\documentclass[a4paper,12pt]{amsart}

\usepackage{amsthm, amsmath, latexsym, amsfonts}
\usepackage[all]{xy}
\usepackage{amssymb}
\usepackage{mathrsfs}

\newtheorem{Lemma}{Lemma}[section]
\newtheorem{Proposition}[Lemma]{Proposition}
\newtheorem{Theorem}[Lemma]{Theorem}

\theoremstyle{remark}
\newcommand{\monbar}{\overline{M}_{0,n}}

\newcommand{\Ker}{\mathrm{Ker}}
\newcommand{\Cox}{\mathrm{Cox}}

\title{A remark on the Cox ring of $\monbar$}

\author{Claudio Fontanari}

\email{fontanar@science.unitn.it}\curraddr{{\sc Dipartimento di Matematica \\  Universit\`a degli Studi di Trento\\ Via Sommarive 14 \\ 38123 Trento \\ Italy.}}

\thanks{This research was partially supported by PRIN 2012 "Geometria delle variet\`a algebriche",
by FIRB 2012 "Moduli spaces and Applications", and by GNSAGA of INdAM (Italy).}



\begin{document}

\begin{abstract}
We present an elementary inductive argument proving that a certain subring of the Cox ring of the moduli space 
$\monbar$ of stable rational curves with $n$ marked points is finitely generated for every $n \ge 3$.
\end{abstract}

\maketitle

\section{Introduction}

Let $\monbar$ be the moduli space of stable rational curves with $n$ ordered distinct marked points. 
After \cite{HuKe}, Question 3.2, it has been an intriguing open problem to establish whether $\monbar$ 
is a Mori dream space. Indeed, an affirmative answer is only known for $n \le 6$ (see \cite{Cas}), 
which is just the range of $n$ for which $\monbar$ is a (smooth) log-Fano variety. On the other hand, 
according to the very recent e-print \cite{CasTev}, $\monbar$ is not a Mori dream space for $n>133$
(at least in characteristic zero). 

Here instead we present an elementary proof of the following fact:

\begin{Theorem}\label{main}   
The graded algebra 
\begin{equation}\label{algebra}
\bigoplus_{c_A \ge 0, c_i \ge 0}  H^0 \left(\monbar, \mathcal{O}_{\monbar}(\sum_{4 \le \vert  A \vert \le (n+1)/2} 
c_A D_A -\sum_i c_i \psi_i) \right)
\end{equation}
is finitely generated for every $n \ge 3$. 
\end{Theorem}

As kindly pointed out to us by Maksym Fedorchuk, such a statement cannot hold 
without any restriction on the index $A$. Indeed, since for every 
$m_A \in \mathbb{Z}$ and for any $m >> 0$ we have 
\begin{eqnarray*}
\sum_A m_A D_A &=& \sum_A m_A D_A + m \sum_j j (n-j) B_j - m(n-1) \sum_i \psi_i = \\
&=& \sum_A c_A D_A -\sum_i c_i \psi_i
\end{eqnarray*}
with all $c_A \ge 0$ and $c_i \ge 0$ (see for instance \cite{Moon}, Lemma 2.9 (2)), it turns out that 
$$
\Cox(\monbar) = \bigoplus_{c_A \ge 0, c_i \ge 0}  H^0 \left(\monbar, \mathcal{O}_{\monbar}(\sum_A c_A D_A -\sum_i c_i \psi_i) \right).
$$

Our inductive argument works for $n \ge 5$ and relies on a vanishing  
$$
H^1 \left( \monbar, \mathcal{O}_{\monbar} (\sum_A c_A D_A - \sum_i c_i \psi_i) \right) = 0
$$
(see Lemma \ref{vanishing} below for a more precise statement). 

We are grateful to Edoardo Ballico, James McKernan and Jenia Tevelev for carefully reading 
a previous version of our manuscript and pointing out to us that such a vanishing definitely 
does not hold for $n=4$, since $\dim H^1(\mathbb{P}^1, \mathcal{O}_{\mathbb{P}^1}(-k))= k-1 > 0$
for $k \ge 2$. We are also grateful to Marco Andreatta, Luca Benzo and Gilberto Bini 
for inspiring conversations on the subject of the present note. 

We work over an algebraically closed field $\mathbb{K}$ of arbitrary characteristic.

\section{The proofs}

The boundary of $\monbar$ is the union of irreducible divisorial components 
\begin{equation}\label{inductive}
D_A = \eta (\overline{M}_{0, A \cup \{ p \} }
\times \overline{M}_{0, P \setminus A \cup \{ q \} }),
\end{equation}
where $P = \{ 1, \ldots, n \}$, $A \subset P$ and 
$\eta $ is the natural morphism obtained by glueing the points $p$ and $q$ into an ordinary node. According to 
\cite{ArbCor}, Lemma 3.3 (see also \cite{Moon}, Lemma 2.5), we have: 

\begin{equation} \label{psi}
\psi_{\vert D_A} = \eta^*(\psi_i) = \left\{
\begin{array}{ll}
(\psi_i, 0) & \mbox{if $i \in A$}\\
(0, \psi_i) & \mbox{if $i \in P \setminus A$}\\ 
\end{array}
\right.
\end{equation}

\noindent
and

\begin{eqnarray} 
\label{self} {D_A}_{\vert D_A} &=& \eta^*(D_A)  = (- \psi_p, -\psi_q) \\ 
\nonumber {D_{P \setminus A}}_{\vert D_A} &=& \eta^*(D_{P \setminus A}) = (- \psi_p, -\psi_q) 
\end{eqnarray}

\noindent
while for every $B$ different from both $A$ and $P \setminus A$:

\begin{equation} \label{boundary}
{D_B}_{\vert D_A} = \eta^*(D_B) = \left\{
\begin{array}{ll}
(D_B, 0) & \mbox{if $B \subset A$}\\
(D_{P \setminus B}, 0) & \mbox{if $P \setminus B \subset A$}\\
(0, D_B) & \mbox{if $B \subset P \setminus A$}\\
(0, D_{P \setminus B}) & \mbox{if $P \setminus B \subset P \setminus A$}\\
(0,0) & \mbox{otherwise.}
\end{array}
\right.
\end{equation}

\begin{Lemma}\label{vanishing}
Let $n \ge 5$. For every $c_A \ge 0$ and every $c_i \ge 0$ we have 
$$
H^1 \left( \monbar, \mathcal{O}_{\monbar} ( \sum_{4 \le \vert  A \vert \le (n+1)/2} c_A D_A - \sum_i c_i \psi_i ) \right) = 0.
$$
\end{Lemma}

\proof 
We argue by induction on $\sum_A c_A$. 

If $\sum_A c_A = 0$, then we need to check that 
\begin{equation}\label{basis}
H^1 \left( \monbar, \mathcal{O}_{\monbar} ( - \sum_i c_i \psi_i ) \right) = 0.
\end{equation}

In order to do so, we make induction on $\sum_i c_i$. Indeed, if $\sum_i c_i = 0$, then 
we have $H^1 ( \monbar, \mathcal{O}_{\monbar} ) =  0$ since $\monbar$ is a smooth rational 
variety and the irregularity is a birational invariant.
Assume now $\sum_i c_i > 0$, so that $c_z > 0$ for some $z$. We know that $\psi_z$ is the class 
of an effective divisor $E$ (see for instance \cite{ArbCor}, (3.7)): 
$$
E = \sum_A D_A
$$ 
with $z \in A$ and $x, y \notin A$ for any choice of distinct elements 
$x, y \in \{ 1, \ldots, n \}$. In particular, $E$ has dimension $\dim(E)=n-4 \ge 1$,
is both reduced and connected (indeed, all such $D_A$ intersect $D_{\{x,y \}}$) and by (\ref{psi}) 
${\psi_z}_{\vert D_A} = (\psi_z, 0)$. Hence from the standard short exact sequence
$$
0 \to \mathcal{O}_{\monbar}(-E) \to \mathcal{O}_{\monbar} \to \mathcal{O}_E \to 0
$$
we deduce 
\begin{eqnarray*}
0 \to \mathcal{O}_{\monbar}(-\sum_i c_i \psi_i) \to \mathcal{O}_{\monbar}
(-\sum_{i \ne z} c_i \psi_i - (c_z -1) \psi_z) \to \\
\to \mathcal{O}_E(-\sum_{i \ne z} c_i \psi_i - (c_z -1) \psi_z) \to 0
\end{eqnarray*}
with 
\begin{eqnarray*}
H^0 \left(\monbar, \mathcal{O}_{\monbar}(-\sum_{i \ne z} c_i \psi_i - (c_z -1) \psi_z) \right)= \\
= H^0 \left(E, \mathcal{O}_E(-\sum_{i \ne z} c_i \psi_i - (c_z -1) \psi_z) \right)
\end{eqnarray*}
both equal to either $\mathbb{K}$ if $\sum_i c_i-1=0$ or to $0$ if $\sum_i c_i-1>0$.
As a consequence, we obtain an exact sequence
\begin{eqnarray*}
0 \to H^1 \left(\monbar, \mathcal{O}_{\monbar}(-\sum_i c_i \psi_i) \right) \to \\
\to H^1 \left(\monbar, \mathcal{O}_{\monbar}(-\sum_{i \ne z} c_i \psi_i - (c_z -1) \psi_z \right)
\end{eqnarray*}
and (\ref{basis}) follows by induction on $\sum_i c_i$.

Assume now $\sum_A c_A > 0$. Among all boundary divisors $D_A$ with $4 \le \vert  A \vert \le (n+1)/2$ 
such that $c_A > 0$, choose $D_B$ minimizing $\vert B \vert$. We may write the short exact sequence
\begin{eqnarray}
\label{short} 0 \to \mathcal{O}_{\monbar}(\sum_{A \ne B} c_A D_A + (c_B-1)D_B -\sum_i c_i \psi_i) \to \\ 
\nonumber \to \mathcal{O}_{\monbar}(\sum_A c_A D_A -\sum_i c_i \psi_i) \to
\mathcal{O}_{D_B}(\sum_A c_A D_A - \sum_i c_i \psi_i) \to 0 
\end{eqnarray}
and deduce
\begin{eqnarray*}
& & H^1 \left(\monbar, \mathcal{O}_{\monbar}(\sum_{A \ne B} c_A D_A + (c_B-1)D_B -\sum_i c_i \psi_i) \right) \to \\ 
&\to& H^1 \left(\monbar, \mathcal{O}_{\monbar}(\sum_A c_A D_A -\sum_i c_i \psi_i) \right) \to \\
&\to& H^1 \left(D_B, \mathcal{O}_{D_B}(\sum_A c_A D_A - \sum_i c_i \psi_i) \right).
\end{eqnarray*}
Now, by inductive assumption on $\sum_A c_A$ we have 
$$
H^1 \left(\monbar, \mathcal{O}_{\monbar}(\sum_{A \ne B} c_A D_A + (c_B-1)D_B -\sum_i c_i \psi_i) \right) = 0.
$$
On the other hand, according to (\ref{inductive}), we have 
\begin{eqnarray*}
& & H^1 \left(D_B, \mathcal{O}_{D_B}(\sum_A c_A D_A - \sum_i c_i \psi_i) \right) = \\
&=& H^1 \left(\overline{M}_{0, B \cup \{ p \} }
\times \overline{M}_{0, P \setminus B \cup \{ q \} }, \eta^*(\sum_A c_A D_A - \sum_i c_i \psi_i) \right)
\end{eqnarray*}
where 
\begin{eqnarray*}
& & \eta^*(\sum_A c_A D_A - \sum_i c_i \psi_i) = \\ 
&=& \left( -c_B \psi_p - \sum_i c_i \psi_i, -c_B \psi_q + \sum_A c_A D_A - \sum_i c_i \psi_i \right)
\end{eqnarray*}
by applying (\ref{psi}), (\ref{self}) and (\ref{boundary}) and recalling our choice of $D_B$ 
such that $B$ minimizes $\vert B \vert$. Hence by K\"unneth formula we obtain 
\begin{eqnarray*}
& & H^1 \left(D_B, \mathcal{O}_{D_B}(\sum_A c_A D_A - \sum_i c_i \psi_i) \right) = \\
&=& H^1 \left(\overline{M}_{0, B \cup \{ p \} }, -c_B \psi_p - \sum_i c_i \psi_i \right) 
\otimes \\
& &\otimes H^0 \left( \overline{M}_{0, P \setminus B \cup \{ q \} }, -c_B \psi_q + \sum_A c_A D_A - \sum_i c_i \psi_i 
\right) \oplus \\
& & \oplus 
H^0 \left(\overline{M}_{0, B \cup \{ p \} }, -c_B \psi_p - \sum_i c_i \psi_i \right) 
\otimes \\
& & \otimes H^1 \left( \overline{M}_{0, P \setminus B \cup \{ q \} }, -c_B \psi_q + \sum_A c_A D_A - \sum_i c_i \psi_i 
\right)
\end{eqnarray*}
with 
$$
H^1 \left(\overline{M}_{0, B \cup \{ p \} }, -c_B \psi_p - \sum_i c_i \psi_i \right)=0
$$
by (\ref{basis}) and 
$$
H^0 \left(\overline{M}_{0, B \cup \{ p \} }, -c_B \psi_p - \sum_i c_i \psi_i \right)=0, 
$$
since $c_B > 0$, hence also
$$
H^1 \left(D_B, \mathcal{O}_{D_B}(\sum_A c_A D_A - \sum_i c_i \psi_i) \right) = 0
$$
and the claim follows.

\endproof

As a consequence of Lemma \ref{vanishing}, from (\ref{short}) we obtain a short exact sequence
\begin{eqnarray*}
0 &\to& H^0 \left(\monbar, \mathcal{O}_{\monbar}(\sum_{A \ne B} c_A D_A + (c_B-1)D_B -\sum_i c_i \psi_i) \right) \to \\ 
&\to& H^0 \left(\monbar, \mathcal{O}_{\monbar}(\sum_A c_A D_A -\sum_i c_i \psi_i) \right) \to \\
&\to& H^0 \left(D_B, \mathcal{O}_{D_B}(\sum_A c_A D_A -\sum_i c_i \psi_i) \right) \to 0
\end{eqnarray*}
for every boundary divisor $D_B$ with $4 \le \vert  B \vert \le (n+1)/2$, so that all restriction maps 
\begin{eqnarray}
\label{restriction} \rho_B: H^0 \left(\monbar, \mathcal{O}_{\monbar}(\sum_{4 \le \vert  A \vert \le (n+1)/2} c_A D_A 
-\sum_i c_i \psi_i) \right) \to \\  
\nonumber \to H^0 \left(D_B, \mathcal{O}_{D_B}(\sum_{4 \le \vert  A \vert \le (n+1)/2} c_A D_A -\sum_i c_i \psi_i) \right) 
\end{eqnarray}
are surjective with kernel
\begin{equation}\label{kernel}
\Ker(\rho_B) = H^0 \left(\monbar, \mathcal{O}_{\monbar}(\sum_{A \ne B} c_A D_A + (c_B-1)D_B -\sum_i c_i \psi_i) \right).
\end{equation}

The idea of the next crucial technical result has been inspired by \cite{Fuj}, Chapter I., Theorem (2.3).  

\begin{Proposition}\label{generation}
If for every $B \subset \{1, \ldots, n \}$ such that $4 \le \vert  B \vert \le (n+1)/2$ the graded algebras 
$$
\mathcal{A}_B = \bigoplus_{c_A \ge 0, c_i \ge 0} H^0 \left(D_B, \mathcal{O}_{D_B}(\sum_{4 \le \vert  A \vert \le (n+1)/2} c_A D_A -\sum_i c_i \psi_i) \right)
$$
are finitely generated by homogeneous elements $\rho_B(\gamma^B_j)$ for $j =1, \ldots,$ $N(B)$, then the graded algebra
$$
\mathcal{A} = \bigoplus_{c_A \ge 0, c_i \ge 0}  H^0 \left(\monbar, \mathcal{O}_{\monbar}(\sum_{4 \le \vert  A \vert \le (n+1)/2} 
c_A D_A -\sum_i c_i \psi_i) \right)
$$   
is finitely generated by (the sections corresponding to) the boundary divisors $D_B$ together with 
the elements $\gamma^B_j$ for every $B$ and $j$.
\end{Proposition}

\proof Let $\mathcal{S} \subseteq \mathcal{A}$ be the subalgebra generated by all sections $s_B$ with 
$D_B = \{ s_B=0 \}$ and by all $\gamma^B_j$. 
In order to prove that $\mathcal{S} = \mathcal{A}$, we are going to check that the following equality holds 
for all homogeneous components: 
\begin{eqnarray*}
\mathcal{S} \cap H^0 \left(\monbar, \mathcal{O}_{\monbar}(\sum_A c_A D_A -\sum_i c_i \psi_i) \right) = \\
= H^0 \left(\monbar, \mathcal{O}_{\monbar}(\sum_A c_A D_A -\sum_i c_i \psi_i) \right).
\end{eqnarray*}

We argue by induction on $\sum_A c_A$. 

\noindent
If $\sum_A c_A = 0$, then 
$H^0 \left(\monbar, \mathcal{O}_{\monbar}(-\sum_i c_i \psi_i) \right)$ equals either $\mathbb{K}$ if 
$\sum_i c_i = 0$ or $0$ if $\sum_i c_i > 0$ and in both cases there is nothing to prove. Thus assume 
$\sum_A c_A > 0$, in particular there is $B$ with $4 \le \vert  B \vert \le (n+1)/2$ such that $c_B > 0$. 
By (\ref{restriction}) and (\ref{kernel}), the vector space 
$$H^0 \left(\monbar, \mathcal{O}_{\monbar}(\sum_A c_A D_A -\sum_i c_i \psi_i) \right)$$ 
is generated by 
$$s_B \otimes H^0 \left(\monbar, \mathcal{O}_{\monbar}(\sum_{A \ne B} c_A D_A + (c_B-1)D_B -\sum_i c_i \psi_i) \right)$$ 
and by a set of elements $\gamma^B_j$ such that $\rho_B(\gamma^B_j)$ generate 
$$H^0 \left(D_B, \mathcal{O}_{D_B}(\sum_A c_A D_A -\sum_i c_i \psi_i) \right),$$ hence 
the claim follows by induction on $\sum_A c_A$. 

\qed

\noindent
\emph{Proof of Theorem \ref{main}.} We argue by induction on $n$, 
starting from the trivial cases $n \le 6$ and then applying 
Proposition \ref{generation} for $n \ge 7$.
 
\endproof

\end{document}